\documentclass[letterpaper, 10 pt, conference]{ieeeconf}  
\IEEEoverridecommandlockouts                              
\usepackage{times} 
\usepackage{subcaption}
\usepackage{amssymb}
\usepackage{bbm}
\let\labelindent\relax         
\usepackage{enumitem}
\usepackage{cite}
\usepackage{tikz}
\usepackage{tkz-graph}
\usepackage{graphicx}  
\usepackage{array} 
\usepackage{multirow}
\usepackage{hhline} 

\usepackage{mathptmx}
\usepackage{amsmath}
\usepackage{mathtools}

\DeclareMathOperator{\im}{im}

\let\leq\leqslant
\let\geq\geqslant
\let\emptyset\varnothing

\newcommand{\calA}{\ensuremath{\mathcal{A}}}

\newcommand{\calC}{\ensuremath{\mathcal{C}}}
\newcommand{\calD}{\ensuremath{\mathcal{D}}}

\newcommand{\calI}{\ensuremath{\mathcal{I}}}

\newcommand{\calL}{\ensuremath{\mathcal{L}}}
\newcommand{\calM}{\ensuremath{\mathcal{M}}}
\newcommand{\calN}{\ensuremath{\mathcal{N}}}
\newcommand{\calO}{\ensuremath{\mathcal{O}}}
\newcommand{\calP}{\ensuremath{\mathcal{P}}}

\newcommand{\calR}{\ensuremath{\mathcal{R}}}
\newcommand{\calS}{\ensuremath{\mathcal{S}}}

\newcommand{\calV}{\ensuremath{\mathcal{V}}}

\newcommand{\hatx}{\ensuremath{\hat{x}}}

\newcommand{\bbR}{\ensuremath{\mathbb{R}}}

\newcommand{\bmat}{\begin{matrix}}
\newcommand{\emat}{\end{matrix}}
\newcommand{\bbm}{\begin{bmatrix}}
\newcommand{\ebm}{\end{bmatrix}}
\newcommand{\bbma}{\begin{bmatrix*}[r]}
\newcommand{\ebma}{\end{bmatrix*}}
\newcommand{\bpm}{\begin{pmatrix}}
\newcommand{\epm}{\end{pmatrix}}
\newcommand{\bvm}{\begin{vmatrix}}
\newcommand{\evm}{\end{vmatrix}}
\newcommand{\bse}{\begin{subequations}}
\newcommand{\ese}{\end{subequations}}
\newcommand{\beq}{\begin{equation}}
\newcommand{\eeq}{\end{equation}}
\newcommand{\beqn}{\begin{equation*}}
\newcommand{\eeqn}{\end{equation*}}

\newcommand{\ben}{\renewcommand{\labelenumi}{(\arabic{enumi})}
\renewcommand{\theenumi}{\arabic{enumi}}\begin{enumerate}}

\newcommand{\een}{\end{enumerate}}

\newcommand{\beni}{\renewcommand{\labelenumi}{(\roman{enumi})}
\renewcommand{\theenumi}{\roman{enumi}}\begin{enumerate}}

\newcommand{\eeni}{\end{enumerate}}

\newcommand{\bena}{\renewcommand{\labelenumi}{(\alph{enumi})}
\renewcommand{\theenumi}{\alph{enumi}}\begin{enumerate}}

\newcommand{\eena}{\end{enumerate}}
\newcommand{\bit}{\begin{itemize}}
\newcommand{\eit}{\end{itemize}}

\newtheorem{thm}{Theorem}
\newtheorem{defn}[thm]{Definition}
\newtheorem{lem}[thm]{Lemma}

\newtheorem{cor}[thm]{Corollary}
\newtheorem{ex}{Example}

\newtheorem{remark}[thm]{Remark}
\newtheorem{prp}[thm]{Proposition}
\newtheorem{pro}{Problem}

\newcommand{\bthm}{\begin{thm}}
\newcommand{\ethm}{\end{thm}}
\newcommand{\blem}{\begin{lem}}
\newcommand{\elem}{\end{lem}}
\newcommand{\bprop}{\begin{prp}}
\newcommand{\eprop}{\end{prp}}
\newcommand{\bas}{\begin{assumption}}
\newcommand{\eas}{\end{assumption}}
\newcommand{\bre}{\begin{remark}}
\newcommand{\ere}{\end{remark}}
\newcommand{\bcor}{\begin{cor}}
\newcommand{\ecor}{\end{cor}}
\newcommand{\bdfn}{\begin{defn}}
\newcommand{\edfn}{\end{defn}}
\newcommand{\bcon}{\begin{conjecture}}
\newcommand{\econ}{\end{conjecture}}
\newcommand{\bex}{\begin{ex}}
\newcommand{\eex}{\end{ex}}
\newcommand{\BP}{\proof}
\newcommand{\EP}{\endproof}
\newcommand{\bpro}{\begin{pro}}
	\newcommand{\epro}{\end{pro}}

\newcommand{\nset}[1]{\ensuremath{\{1,2,\ldots,#1\}}}

\title{\LARGE \bf 
	Fault detection and isolation for linear structured systems}

\author{Jiajia~Jia,
	Harry~L.~Trentelman~
	and~M.~Kanat~Camlibel 
	\thanks{The authors are with the Bernoulli Institute for Mathematics, Computer Science and Artificial Intelligence, University of Groningen,  9700 Groningen, The Netherlands. {\tt\small (e-mail: \{j.jia, h.l.trentelman, m.k.camlibel\}@rug.nl.)}}
}
\begin{document}
	
	\maketitle
	\thispagestyle{empty}
	\pagestyle{empty}

\begin{abstract}                
This paper deals with the fault detection and isolation (FDI) problem for linear structured systems in which the system matrices are given by zero/nonzero/arbitrary pattern matrices.
In this paper, we follow a geometric approach to verify solvability of the FDI problem for such systems.
To do so, we first develop a necessary and sufficient condition under which the FDI problem for a given particular linear time-invariant system is solvable.
Next, we establish a necessary condition for solvability of the FDI problem for linear structured systems.
In addition, we develop a sufficient algebraic condition for solvability of the FDI problem in terms of a rank test on an associated pattern matrix.
To illustrate that this condition is not necessary, we provide a counterexample in which the FDI problem is solvable while the condition is not satisfied.
Finally, we develop a graph-theoretic condition for the full rank property of a given pattern matrix, which leads to a graph-theoretic condition for solvability of the FDI problem.
\end{abstract}

\section{Introduction}\label{s:int} 
This paper is concerned with the FDI problem for  linear time-invariant (LTI) systems with faults. 
This problem has received considerable attention within the control community in the past decades and this has lead to several approaches to FDI, see, e.g., \cite{C1980, M1986, MVW1989, FRANK1996, CD2007, REC2015} and the references therein.
Among these references, those closer to the results presented in the current paper are \cite{M1986} and \cite{REC2015}, in which FDI for LTI systems is performed using unknown input observers that enable so-called output separability of the fault subspaces.
If such observers exist, then we say that for the given system the FDI problem is solvable.

Although conditions for solvability of the FDI problem for a given LTI system have been introduced in \cite{M1986}, their application relies on the exact knowledge
of the dynamics of this system, meaning that precise information on the system matrices is required. However, in many scenarios, such knowledge is unavailable, and only the zero/nonzero/arbitrary structure can be acquired.
This leads to the concept of linear structured system introduced in \cite{JWTC2019} which represents a family of LTI systems sharing the same structure.
A large amount of literature has been devoted to analyzing system-theoretical properties for linear structured systems.
For instance, strong structural controllability has been studied in \cite{CM2013, MZC2014, TD2015,JWTC2019}, strong targeted controllability in \cite{MCT2015, vWCT2017}, and identifiability in\cite{WTC2018}.

Roughly speaking, in the framework of linear structured systems, the research on the FDI problem can be subdivided into two directions.
The first direction aims at providing conditions under which the FDI problem is solvable for at least one member of a given structured system, see, e.g.,  \cite{CD2007, CDSM2000, CDA2008}.
The other direction aims at establishing conditions to guarantee that the FDI problem is solvable for all members of a given structured system, see, e.g., \cite{REC2015}. 
In the present paper, we will pursue the second research direction.
For a given structured system, if the FDI problem for all systems in the structured system is solvable, then we say that the FDI problem for this structured system is {\em solvable}.
To the best of our knowledge,  in this direction the only existing work is \cite{REC2015}, which has studied a special kind of linear structured system, named systems defined on graphs.
The goal of the present paper is to provide conditions under which the FDI problem is solvable for a general structured system. 
The main contributions of this paper are the following:
\begin{enumerate}
	\item We develop a necessary and sufficient condition under which the FDI problem is solvable for a given particular LTI system.               
	\item Based on the condition for a particular LTI system, we establish a necessary condition for solvability of the FDI problem for a given structured system.
	         Next, we develop a  sufficient algebraic condition.
              This condition is given in terms of a rank test on a pattern matrix  associated with the structured system.
              Moreover, we provide a counterexample to show that this condition is not necessary.
	\item Using the concept of  colorability of a graph, we provide a graph-theoretic condition for solvability of the FDI problem for a given structured system. 
\end{enumerate}

This paper is structured as follows. 
In Section \ref{s:pre}, we review concepts and preliminary results on geometric control theory and the geometric approach to the FDI problem for particular LTI systems. 
In addition, we introduce the concept of linear structured systems and formulate the problem studied in this paper.  
Section \ref{s:single} presents a necessary and sufficient condition under which for a given particular LTI system the FDI problem is solvable.
Section \ref{s:structure} provides a necessary and a sufficient algebraic condition for solvability of the FDI problem for structured systems.
Next, in Section \ref{s:graph} we  establish a graph-theoretic condition for solvability of the FDI problem.
Finally, section \ref{s:conclusions} concludes this paper.

\section{Preliminaries and problem statement}\label{s:pre}
Let $\mathbb{R}$ and $\mathbb{R}^{n}$ denote the field of real numbers and the vector space of $n$-dimensional real vectors, respectively. 
Likewise, we denote the space of $n \times m$ real matrices by $\mathbb{R}^{n \times m}$.
For a given matrix $M \in \mathbb{R}^{n \times m}$, the $i$th column of $M$ is denoted by $M_i$.
Moreover, $I$ and $0$ will denote  identity and zero matrices of appropriate dimensions, respectively.

\subsection{Geometric control theory}\label{s:gct}
Geometric control theory plays a fundamental role in this paper. 
Therefore, in this subsection, we will give a brief review of some basic concepts in this field.   
Consider the LTI system
\beq \label{eq:LTI}
	 \begin{split}
		\dot{x} &=  Ax + Bu \\ 
		y &=  Cx,
	\end{split}
\eeq
where $x \in \bbR^n$, $u \in \bbR^m$ and $y \in \bbR^p$ are the state, input and output, respectively, and
$A$, $B$ and $C$ are matrices of appropriate dimensions.
A subspace $\calS \subseteq \bbR^n$ is called {\em $(C,A)$-invariant} if 
$A (\calS \cap \ker C) \subseteq \calS.$
This condition is equivalent to the existence of a matrix $G \in \bbR^{n \times p}$ such that $(A + GC)\calS \subseteq \calS.$
Such a $G$ is called a \emph{friend} of $\calS$. 
A family $ \{\calS_i\}^{k}_{i = 1} $ of $(C,A)$-invariant subspaces of $\bbR^n$ is called \emph{compatible} if the subspaces $\calS_i$ have a common friend.
Given the system \eqref{eq:LTI}, a family of subspaces $\{\calS_i \}^{k}_{i = 1}$ is called \emph{output separable}   if for $i = 1,2, \ldots,k$
\beqn \label{eq:outputseparable}
C\calS_i \cap (\sum_{ j \neq i} C \calS_j) = \{ 0 \}.
\eeqn
Any output separable family of $(C,A)$-invariant subspaces is compatible \cite[Lemma 2]{M1986}.
Moreover, if it also satisfies the condition that $C\calS_i \neq \{0\}$ for $i = 1, 2, \ldots, k$, we say that the family $\{C\calS_i \}^{k}_{i = 1}$ is {\em independent}.

For a given subspace $\calD \subseteq \bbR^n$, there exists a smallest $(C,A)$-invariant subspace containing $\calD$, denoted by $\calS^{\ast}$.
Such a minimal subspace can be computed by the following subspace algorithm (see, e.g., the conditioned invariant subspace algorithm p.111 of \cite{TSH2012}):
\beq \label{eq:algorithm1}
\begin{split}
	\calS^0 &=  \calD \\ 
	\calS^k &=  \calD + A(\calS^{k-1} \cap \ker C) \mbox{ for } k = 1, 2, \ldots.
\end{split}
\eeq
Denote the dimension of $\calD$ by $\dim (\calD)$.
It follows from Theorem 5.8 of \cite{TSH2012} that there exists $k \leq n - \dim \calD$ such that $\calS_{k} = \calS_{k+1}$, and hence $\calS^{\ast} = \calS_{k}$.

\subsection{The geometric approach to the FDI problem for LTI systems}\label{s:GFDI}

In this subsection, we will review the geometric approach to the FDI problem for LTI systems.
Consider the LTI system 
\beq \label{eq:LTIf}
\begin{split}
    \dot{x} &=  Ax + Lf \\ 
    y &=  Cx,
\end{split}
\eeq
where $x \in \bbR^n$,  $f \in \bbR^q$ and $y \in \bbR^p$ are the state,  fault and output, respectively, and $A$,   $L$ and $C$ are matrices of appropriate dimensions.
We denote the system \eqref{eq:LTIf} by $(A,L,C)$. 
We say that the $i$th fault occurs if $f_i \neq 0$ (i.e., not identically equal to $0$), where $f_i$ is the $i$th component of $f$.
Following the approach proposed in \cite{M1986}, the FDI problem for \eqref{eq:LTIf} amounts to finding $G \in \bbR^{n \times p}$
such that the family of subspaces $\{C\calV_i\}_{i = 1}^{q}$ is independent, where $\calV_i$ is the smallest $(A+GC)$-invariant subspace containing $\im L_i$.
If such $G$ exists, then we say that the FDI problem is {\em solvable}.
In what follows, we will  briefly explain this approach. 
Suppose that we have found a $G$ satisfying the above constraints. 
Consider the state observer
\beq \label{eq:observer}
\dot{\hatx}  = (A + GC)\hatx - Gy.
\eeq
Define the {\em innovation} as 
\[
r \coloneqq C\hatx - y
\]
and 
{\em error}  
\[
e \coloneqq \hatx -x.
\]
By interconnecting \eqref{eq:LTIf} and \eqref{eq:observer}, we obtain
\begin{equation} \label{eq:ES}
\begin{split}
	\dot{e} & = (A + GC)e - Lf\\ 
    r & =   Ce.
\end{split}
\end{equation}
Note that in this paper, we do not consider any stability requirement on the observer, which means that we do not require $e(t) \rightarrow 0$, and  we assume that $e(0) = 0$.
Under this assumption, for any fault $f$, the resulting error trajectory $e(t)$ lies in the reachable subspace of $(A+GC,L)$, which is clearly equal to 
$\calV_1+ \calV_2 + \cdots+\calV_q$.
For the corresponding innovation trajectory $r(t)$ we then have $$r(t) \in C\calV_1+ C\calV_2 +\cdots +C\calV_q.$$
If the family $ \{C\calV_i\}^q_{i=1}$ is independent, then this is a direct sum, and $r(t)$ can be written uniquely as 
\beq\label{eq:rt}
r(t) = r_1(t) + r_2(t) + \cdots +r_q(t)
\eeq
with  $r_i(t) \in C\calV_i$ for all $t$.
The unique representation \eqref{eq:rt} can be used to determine whether the $i$th fault occurs.
Indeed in  \eqref{eq:rt} $r_i \neq 0$ (i.e., not identically equal to $0$) only if  $f_i \neq 0$.
To see this, note that $f_i(t) = 0$ for all $t$ implies $e(t) \in \sum_{ j \neq i} \calV_j$, so $r(t) \in \sum_{ j \neq i} C\calV_j$, equivalently, $r_i(t) = 0$ for all $t$.

Let  $\calS_i^\ast$ be the smallest $(C,A)$-invariant subspace containing $\im L_i$.
In \cite{M1986} it has been shown that the FDI problem for the system \eqref{eq:LTIf} is solvable if and only if the family $\{C\calS^\ast_i\}_{i = 1}^q$ is independent, i.e.,
the family $\{\calS^{\ast}_i\}_{i = 1}^{q}$ is output separable and $C\calS^\ast_i \neq \{0\}$ for $i = 1, 2, \ldots, q.$

\subsection{Linear structured systems and problem formulation}\label{s:problem}
Again, consider the LTI system \eqref{eq:LTIf}.
In many scenarios, the exact values of the entries in the system matrices are not known, but some entries are known to be always zero, some are nonzero, and the remaining entries are arbitrary real numbers.
To describe such kind of matrices, the authors in \cite{JWTC2019} have introduced the definition of {\em pattern matrix} as follows.

A pattern matrix is a matrix with entries in the set of symbols $\{0,\ast,?\}$.
The set of all $r \times s$ pattern matrices is denoted by $\{0,\ast,?\}^{r \times s}$.
For a given $r \times s$ pattern matrix $\mathcal{M}$, we define the \emph{pattern class} of $\mathcal{M}$ as
\begin{equation*}
\begin{aligned}
\mathcal{P}(\mathcal{M}) := \{M \in \mathbb{R}^{r \times s} \mid  & M_{ij} = 0 \text{ if } \mathcal{M}_{ij} = 0, \\& M_{ij} \neq 0 \text{ if } \mathcal{M}_{ij} = \ast \}.
\end{aligned}
\end{equation*}
This means that for a matrix $M \in \mathcal{P}(\mathcal{M})$, the entry $M_{ij}$ is either (i) \emph{zero} if $\mathcal{M}_{ij} = 0$, (ii) \emph{nonzero} if $\mathcal{M}_{ij} = \ast$, or (iii) \emph{arbitrary} (zero or nonzero) if $\mathcal{M}_{ij} = \: ?$.

Let $\calA \in \{0,\ast,?\}^{n \times n}$,  $\calL \in \{0,\ast,?\}^{n \times q}$ and $\calC \in \{0,\ast,?\}^{n \times p}$.
The family of systems $(A,L,C)$ with $A \in \calP(\calA)$,  $L \in \calP(\calL)$ and $C \in \calP(\calC)$ is called the {\em linear structured system  associated with $\calA$, $\calL$, and $\calC$}.
Throughout this paper, we use $(\calA,\calL,\calC)$ to represent this structured system, and we write $(A,L,C) \in (\calA,\calL,\calC)$ if $A \in \calP(\calA)$,  $L \in \calP(\calL)$ and $C \in \calP(\calC)$. 
Based on these notions and notations, we define the FDI problem for $(\calA,\calL,\calC)$ to be {\em solvable} if the FDI problem is solvable for every $(A,L,C) \in (\calA,\calL,\calC)$.
The research problem of this paper is then formally stated as follows. 
\begin{pro}
Given $(\calA,\calL,\calC)$, find conditions under which  the FDI problem is solvable for  $(\calA,\calL,\calC)$.
\end{pro}

\section{A necessary and sufficient condition for solvability of the FDI problem for $(A,L,C)$}\label{s:single}

 In this section, we will establish a necessary and sufficient condition under which  the FDI problem is solvable for a given LTI system $(A,L,C)$ of the form \eqref{eq:LTIf}.
Recall that solvability of the FDI problem for $(A,L,C)$ is equivalent to the independence of the family $\{C\calS^\ast_i\}_{i = 1}^q$, where $\calS_i^\ast$ is the smallest $(C,A)$-invariant subspace containing $\im L_i$ ($i = 1, 2, \ldots,q$).
Therefore, we will first provide a characterization of $\calS_i^\ast$.
Let $d_i$ be a positive integer such that 
\[ CA^j L_i = 0 \text{ for } j = 0,1, \ldots, d_i - 2 ~ \text{ and}~~ CA^{d_i-1} L_i \neq 0.\]
Here and in the sequel, we define $A^0 := I$.
It is obvious from the Cayley-Hamilton theorem that either $d_i \leq n$ or $d_i$ does not exist.
If this $d_i$ exists, we then call it the {\em index} of $(A,L_i,C)$.

We are now ready to state a characterization of  $C\calS_i^\ast$ in the following lemma.
\blem \label{l:subspace}
Consider the system $(A,L,C)$ of the form \eqref{eq:LTIf}. Let $i \in \{1,2,\ldots,q\}$.
Denote by $\calS^{\ast}_i$ the smallest $(C,A)$-invariant subspace containing $\im L_i$.
Then, we have that 
\beqn
C\calS^{\ast}_i  =
\begin{cases}
 \im C A^{d_i -1} L_i &  ~~~\text{if the index $d_i$ of $(A, L_i,C)$exists,}\\
 \{0\}& ~~~\text{otherwise.}
\end{cases}
\eeqn
\elem
\vspace{1.9pt}

\proof
In this proof, we will employ the recurrence relation \eqref{eq:algorithm1} to prove the statement.
Let $\calS^\ell_i$ be the sequence of subspaces given by
\beq \label{eq:0}
\begin{aligned}
\calS^0_i &= \im L_i ,\\
\calS^\ell_i &= \im L_i + A(\calS^{\ell -1}_i \cap \ker C) \mbox{ for } \ell = 1,2, \ldots.
\end{aligned}
\eeq
We then distinguish two cases: $(i)$ $d_i$ exists, and $(ii)$ $d_i$ does not exist.

In  case $(i)$, we have that  
 \beq \label{eq:2}
 C A^k L_i = 0 \mbox{ for } k = 0, 1, \ldots, d_i-2
 \eeq 
and
\beq \label{eq:1}
C A^{d_i - 1} L_i \neq 0.
\eeq
By combining \eqref{eq:0} and \eqref{eq:2}, it can be verified directly that
\beq \label{eq:3}
 \calS^k_i = \im \bbm L_i & A L_i & \ldots & A^{k} L_i\ebm \text{ for } k = 0, 1, \ldots, d_i -1.
\eeq
Now, we claim that:
\begin{enumerate}[label={(\alph*)}]
\item $\calS^{d_i -1}_i = \calS^{d_i}_i,$
\item the dimension of $\calS^{d_i -1}_i$ is strictly larger than that of $\calS^{d_i -2}_i$.
\end{enumerate}
If both claims (a) and (b) are true, then $$\calS^\ast_i = \calS^{d_i -1}_i = \im \bbm L_i & A L_i & \ldots & A^{d_i - 1} L_i\ebm,$$ and hence $C\calS^\ast_i = \im C A^{d_i - 1}L_i.$ 
Note that (a) follows immediately from \eqref{eq:1} and \eqref{eq:3}:
$$
\begin{aligned}
\calS^{d_i -1}_i &\stackrel{\eqref{eq:3}}{=}  \im \resizebox{0.38\hsize}{!}{ $\bbm L_i & A L_i & \ldots & A^{d_i - 1} L_i\ebm$}\\
\calS^{d_i}_i &= \im L_i + A(\calS^{d_i -1}_i \cap \ker C)\\ &\stackrel{\eqref{eq:1}}{=}  \im \resizebox{0.38\hsize}{!}{ $\bbm L_i & A L_i & \ldots & A^{d_i - 1} L_i\ebm$}.
\end{aligned}
$$
To prove (b), we assume that (b) is not true, i.e., $$\calS^{d_i -1}_i = \calS^{d_i -2}_i = \im \resizebox{0.38\hsize}{!}{ $\bbm L_i & A L_i & \ldots & A^{d_i - 2} L_i\ebm$}.$$
This implies  $$A^{d_i - 1} L_i \in  \im \bbm L_i & A L_i & \ldots & A^{d_i - 2}  L_i\ebm \subseteq \ker C,$$
which contradicts \eqref{eq:1}, and hence (b) is proved.

For case (ii), we have 
 \beq \label{eq:4}
 C A^k L_i = 0 \mbox{ for } k = 0, 1, \ldots, n-1.
 \eeq 
By combining \eqref{eq:0} and \eqref{eq:4}, we obtain 
\[ \begin{aligned}\calS^{n-1}_i &= \im \bbm L_i & AL_i & \cdots & A^{n-1}L_i \ebm \subseteq \ker C\\
                         \calS^{n}_i &= \im \bbm L_i & AL_i & \cdots & A^{n-1}L_i & A^{n}L_i \ebm\end{aligned}.\]
It then follows from the Caley-Hamilton theorem that $A^{n}L_i \in \calS^{n-1}_i$, i.e., $\calS^{n-1}_i = \calS^{n}_i,$ and hence 
$\calS^\ast_i = \calS^{n-1} \subseteq \ker C.$
Therefore, we have $C\calS^\ast_i = \{0\}.$
This completes the proof.
\endproof
\vspace{1.5pt}

By the above lemma,  the family $\{C\calS^\ast_i\}_{i = 1}^q$ of subspaces is independent if and only if  the index $d_i$ exist for $i = 1,2, \ldots,q$, and 
the vectors $\{CA^{d_i -1} L_i\}_{i = 1}^q$ are linearly independent.
Thus we arrive at the main result of this section which provides a necessary and sufficient condition under which the FDI problem for $(A,L,C)$ is solvable.
\vspace{1.5pt}
\bthm \label{t:ROut}
Consider the system $(A, L, C)$ of the form \eqref{eq:LTIf}. 
The FDI problem for $(A,L,C)$ is solvable if and only if the index $d_i$ exists for $i = 1,2, \ldots,q$, and the matrix $R$ has full column rank, where $R$ is defined by
\beq \label{eq:R}
R \coloneqq  \bbm CA^{d_1 -1}L_1 & CA^{d_2 -1}L_2 & \cdots & CA^{d_q -1}L_q \ebm.
\eeq
\ethm
\vspace{3pt}
\BP
The proof follows immediately from Lemma \ref{l:subspace}  and is hence omitted. 
\EP
\vspace{1.5pt}
\section{Algebraic conditions for solvability of the FDI problem for $(\calA,\calL,\calC)$}\label{s:structure}

In this section, we will establish a necessary condition and a sufficient condition that enables the FDI problem for a given structured system $(\calA, \calL, \calC)$ to be solvable.
Before presenting the results of this section, we first provide some background on operations on pattern matrices. More details can be found in \cite{B2019}.
Addition and multiplication within the set $\{0,\ast,?\}$ are defined in Table \ref{ta:results} below.

\begin{table}[h]
	\caption{Addition and multiplication within the set $\{0,\ast,?\}$.}
		\label{ta:results}
		\begin{center}
	\begin{tabular}{|c|ccc|ccc|}
		\hline
		$+$ & $0$ & $\ast$ & $?$  \\  \hline 
		$0$& $0$ & $\ast$ & $?$  \\ 
		$\ast$& $\ast$ & $?$ & $?$  \\ 
		$?$& $?$ & $?$ & $?$  \\ \hline
	\end{tabular}
	~~
	\begin{tabular}{|c|ccc|ccc|}
		\hline
		$\boldsymbol{\cdot}$  & $0$ & $\ast$ & $?$ \\  \hline 
		$0$&  $0$ & $0$ & $0$ \\ 
		$\ast$& $0$ & $\ast$ & $?$ \\ 
		$?$&  $0$ & $?$ & $?$ \\
		 \hline
	\end{tabular}
	\end{center}
\end{table}
Based on the operations in Table \ref{ta:results}, multiplication of pattern matrices is then defined as follows.
\bdfn \label{eq:mPm}
Let  $\calM \in \{0,\ast,?\}^{r \times s}$ and $\calN \in \{0,\ast,?\}^{s \times t}$.
The product of $\calM$ and $\calN$ is defined as $\calM \calN \in \{0,\ast,?\}^{r \times t}$ given by 
\beq
(\calM \calN)_{ij} :=  \sum_{k = 1}^{q} (\calM_{i k} \boldsymbol{\cdot} \calN_{k j
}) ~~ i = 1,2, \ldots, r, ~~ j = 1,2, \ldots, t.
\eeq
\edfn

It is easily seen that $MN \in \calP(\calM \calN)$ for every pair of matrices $M \in \calP(\calM)$ and $N \in \calP(\calN)$. 
If $ r = s$,  we call $\calM$ a square pattern matrix.
For any given non-negative integer $k$, we define the $k$th power $\calM^k$ recursively by
\[
\calM^0 = \calI, ~~ \calM^i = \calM^{i-1}\calM, ~~i = 1,2, \ldots, k,
\]
where $\mathcal{I}$ represents a square pattern matrix of appropriate dimensions with all diagonal entries equal to $\ast$ and all off-diagonal equal to $0$.
In the sequel, let $\mathcal{O}$ denote any pattern matrix of appropriate dimensions with all entries equal to $0$.

Next, consider the system $(\calA,\calL,\calC)$.
Let $\calL_i$ represent the $i$th column of $\calL$ for $i = 1, 2, \ldots,q$.
Let $\eta_i$ be a positive integer such that 
\[\calC \calA^{j} \calL_i = \calO \text{ for } j = 0,1, \ldots, \eta_i -2 ~\text{ and }~ \calC \calA^{\eta_i -1} \calL_i \neq \calO.\]
If $\eta_i$ exists, then we call it the {\em index} of $(\calA,\calL_i,\calC)$. 
In the sequel, we will write $(A,L_i,C) \in (\calA,\calL_i,\calC)$ if  $A \in \calP(\calA)$,  $L_i \in \calP(\calL_i)$ and  $C \in \calP(\calC)$.
Before continuing to explore conditions for solvability of the FDI problem for $(\calA,\calL,\calC)$, we first provide the following lemma which states  the relationship between the index of $(\calA,\calL_i,\calC)$ and that of $(A,L_i,C) \in (\calA,\calL_i,\calC)$.
\blem \label{l:eta}
Consider the pattern matrix triple $(\calA,\calL_i,\calC)$.
Then the following holds:
\beni
\item Let $ (A, L_i, C) \in (\calA,\calL_i,\calC)$. If both the index $\eta_i$ of $(\calA,\calL_i,\calC)$ and the index $d_i$ of $(A, L_i, C)$ exist, then $d_i \geq \eta_i$.
\item Suppose that the index $\eta_i$ of $(\calA,\calL_i,\calC)$ exists, and suppose further that at least one entry of  $\calC \calA^{\eta_i -1}\calL_i$ is equal to $\ast$. 
Let $(A, L_i, C) \in (\calA,\calL_i,\calC)$.
Then, the index $d_i$ of $(A, L_i, C)$ exists and $d_i = \eta_i$.
\item If the index of $(\calA,\calL_i,\calC)$ does not exist, then  the index of $(A, L_i, C) $ does not exist for any $(A, L_i, C) \in (\calA,\calL_i,\calC)$.
 \eeni
\elem
\BP
By Definition \ref{eq:mPm}, it follows that the vector $CA^{\ell}L_i \in \calP(\calC\calA^{\ell}\calL_i)$ for $i = 0, 
1, \ldots$ and for all $(A,L_i,C) \in (\calA, \calL_i,\calC)$.
In order to prove (i), suppose that both the index $\eta_i$ of $(\calA,\calL_i,\calC)$ and the index $d_i$ of $(A, L_i, C)$ exist.
By the definition of  $\eta_i$ we have that $\calC \calA^{\ell} \calL_i = \calO$ for $\ell = 0, 1, \ldots, \eta_i -2$, and by the definition of $d_i$ it follows that $CA^{d_i-1} L_i \neq 0$. 
Therefore, we obtain $d_i \geq \eta_i$.
Next, to prove (ii), we assume that $\calC \calA^{\eta_i -1}\calL_i$ contains at least one $\ast$ entry, which implies that all the vectors in the pattern class $\calP(\calC \calA^{\eta_i -1} \calL_i)$ are unequal to $0$. 
Let $(A, L_i, C) \in (\calA,\calL_i,\calC)$. Clearly, the vector $CA^{\eta_i -1}L_i \in \calP(\calC \calA^{\eta_i -1} \calL_i)$, and hence $CA^{\eta_i -1} L_i \neq 0$.
By definition, the index $d_i$ of $(A, L_i, C)$ must exist and $d_i \leq \eta_i$.
Recalling (i), we conclude that $d_i = \eta_i$.
The proof of (iii) is trivial. Indeed, suppose that the index of $(\calA,\calL_i,\calC)$ does not exist. It then follows that  $\calC \calA^{\ell} \calL_i = \calO$ for  $\ell = 0 ,1, \ldots$, which implies that $CA^{\ell}L_i$ is equal to $0$ for every  $(A, L_i, C) \in (\calA,\calL_i,\calC)$. 
That is, the index of  $(A, L_i, C) $ does not exist for any $(A, L_i, C) \in (\calA,\calL_i,\calC)$.
\EP

To illustrate the above lemma, we now provide an example.
\bex \label{ex:1}
Consider the system $(\calA, \calL, \calC)$ with
\beq \label{eq:ex1}
\calA = \bbm 0 & 0   & 0 \\
\ast  & 0   & 0 \\ 
0  & 0   & 0 \ebm,~ 
\calL = \bbm \ast & 0 & 0  \\ 
0 &\ast & 0 \\
0 &\ast & \ast\ebm,      
~ 
\calC = \bbm ? & \ast   & 0   \\
0    & \ast  & 0 \ebm.
\eeq
Let $\calL_1$, $\calL_2$ and $\calL_3$ denote the first, second and third column of $\calL$.
For $\calL_1$ and $\calL_2$ we compute
\[\calC \calL_1 = \bbm ? \\ 0\ebm \neq \calO \text{ and } \calC \calL_2 = \bbm \ast \\ \ast\ebm \neq \calO.\]
This implies that $\eta_1 = \eta_2=1$, where $\eta_i$ is the index of $(\calA,\calL_i,\calC)$ for $i = 1,2$.
In addition, for $\calL_3$ we compute
\[
\calC \calA^{\ell}\calL_3 = \calO \text{ for } i = 0, 1, 2, \ldots
\]
which implies that the index of $(\calA,\calL_3,\calC)$ does not exists. 
Next, we will show that for some $(A,L_1,C) \in (\calA,\calL_2,\calC)$ the index $d_1$ of $(A,L_1,C)$ is larger than $\eta_1$, for every $(A,L_2,C) \in (\calA,\calL_2,\calC)$ its index $d_2$ is equal to $\eta_2$, and for every $(A,L_3,C) \in (\calA,\calL_3,\calC)$ its index does not exists,.
Indeed, for $A \in \calP(\calA)$, $L \in \calP(\calL)$ and $C \in \calP(\calC)$ we have
\beq \label{eq:ALC}
A = \resizebox{0.21\hsize}{!}{ $\bbm 0 & 0   & 0 \\
c1  & 0   & 0 \\ 
0  & 0   & 0 \ebm$},~ 
L = \resizebox{0.24\hsize}{!}{ $\bbm c_2 & 0   & 0\\ 
0 & c_3 &0\\
0 & c_4 & c_5 \ebm $},     
~
C = \resizebox{0.24\hsize}{!}{ $\bbm \lambda_1 & c_6   & 0   \\
0  & c_7   & 0 \ebm$},
\eeq
where $c_1, c_2, \ldots, c_7$ are arbitrary nonzero real numbers, and $\lambda_1$ is an arbitrary real number.
Next, we compute
\beq \label{eq:CL}
\bbm CL_1 & CL_2\ebm = \bbm  \lambda_1 c_2 & c_3 c_6 \\ 0 & c_3 c_9\ebm \mbox{ and }~ CAL_1= \bbm c_1 c_2 c_6\\ c_1 c_2 c_7  \ebm.
\eeq
Thus,  for all choices of $c_1, c_2, \ldots, c_7$ and $\lambda_1$ we have  $d_2 = 1 = \eta_2$, while if $\lambda_1  = 0$ then $d_1 = 2 > \eta_1$ and  otherwise  $d_1 = 1 = \eta_1$.
In addition, it is obvious that for all choices of $c_1, c_2, \ldots, c_7$ and $\lambda_1$ we have $C A^\ell L_3 = 0$ for $\ell = 0, 1, \ldots$, and hence the index of $(A,L_3,C)$ does not exist.
\eex

Lemma \ref{l:eta} immediately yields a necessary condition for solvability of the FDI problem for  $(\calA,\calL,\calC)$. 
\bthm \label{t:rankoutputseparable}
Consider the system $(\calA,\calL,\calC)$.
Suppose that the FDI problem for $(\calA, \calL, \calC)$ is solvable.
Then, the index $\eta_i$ of $(\calA,\calL_i,\calC)$ exists for all $i = 1,2,\ldots q$.
\ethm
\BP
Since the FDI problem for $(\calA, \calL, \calC)$ is solvable,  the FDI problem is solvable for all $(A,L,C) \in (\calA,\calL,\calC)$.
Assume that for some $i \in \nset{q}$ the index $\eta_i$ of $(\calA,\calL_i,\calC)$ does not exist.
By statement (iii) of Lemma \ref{l:eta}, it follows that the index $d_i$ of $(A,L_i,C)$ does not exist for any $(A,L_i,C) \in (\calA,\calL_i,\calC).$
It then follows from Theorem \ref{t:ROut} that the FDI problem for  $(A,L,C)$ is not solvable for any $(A,L,C) \in (\calA,\calL,\calC).$
Therefore, we reach a contradiction and complete the proof. 
\EP

By the above theorem, in the sequel we will assume that for all $i = 1,2,\ldots q$ the indices $\eta_i$ exist.
Based on this assumption, we will continue to explore sufficient conditions for solvability of the FDI problem for  $(\calA,\calL,\calC)$.
To do so, we first define the following  pattern matrix associated with  $(\calA, \calL, \calC)$:
\beq \label{eq:calR}
\calR \coloneqq 
\bbm 
\calC \calA^{\eta_1 - 1} \calL_1 & \calC \calA^{\eta_2 - 1} \calL_2 & \cdots & \calC \calA^{\eta_q -1} \calL_q
\ebm,
\eeq
where $\eta_i$ is the index of $(\calA,\calL_i,\calC)$.
We say that $\calR$ has {\em full column rank} if all the matrices in the pattern class $\calP(\calR)$ have full column rank.
We are now ready to establish a sufficient condition for solvability of the FDI problem for  $(\calA,\calL,\calC)$. 
\bthm \label{t:rankoutputseparable}
Consider the system $(\calA,\calL,\calC)$.
Let $\calR$ be the pattern matrix given by \eqref{eq:calR}.
The FDI problem for $(\calA, \calL, \calC)$ is solvable if $\calR$ has full column rank.
\ethm
\begin{proof}
Since $\calR$ has full column rank, each column of $\calR$ contains at least one $\ast$ entry.
Let $(A,L,C) \in (\calA,\calL,\calC).$
By (ii) of Lemma \ref{l:eta} it follows that $d_i = \eta_i$, where $d_i$ is the index of $(A,L_i,C)$ for $i = 1, 2, \ldots, q$.
This implies that the matrix $R$ given by \eqref{eq:R} is in $\calP(\calR)$, and hence $R$ has full column rank.
It then follows from Theorem \ref{t:ROut} that the FDI problem is  solvable.
Since $(A,L,C)$ is an arbitrary system in $(\calA,\calL,\calC)$, we conclude that the FDI problem for $(\calA,\calL,\calC)$ is solvable and complete the proof.
\end{proof}


Note that the condition given in Theorem  \ref{t:rankoutputseparable} is sufficient but not necessary.
To show this, we provide the following counterexample.
\bex \label{ex:2}
Consider the system $(\calA, \calL, \calC)$ with
\[ \calA = \bbm 0 & 0 \\
0  & 0 \ebm, ~~
\calL = \bbm \ast & \ast   \\ 
0 & \ast \ebm       
,~~
\calC = \bbm \ast & \ast \\
\ast  & 0  \ebm.\]
Let $\calL_1$ and $\calL_2$ be the first and second column of $\calL$.
We compute 
\[
\calC \calL_1 = \bbm \ast \\ \ast\ebm \text{ and } \calC \calL_2 = \bbm  ? \\  \ast\ebm,
\]
and, by \eqref{eq:calR}, $\calR = \bbm \ast & ? \\ \ast & \ast \ebm$.
Since $\bbm 1 & 1\\ 1& 1 \ebm \in \calP(\calR)$, $\calR$ does not have full column rank. 
Next, we will show that, however, the FDI problem for $(\calA,\calL,\calC)$ is solvable.
Due to Theorem \ref{t:ROut}, it suffices to show that for each $(A,L,C) \in (\calA,\calL,\calC)$ the associated matrix $R$ has full column rank.
Clearly, every $(A,L,C) \in (\calA,\calL,\calC)$ has the form
\[A = \bbm 
0 & 0\\ 0& 0
\ebm,~~
L = \bbm c_1 & c_2   \\ 
0 & c_3  \ebm  ,~~    
C = \bbm c_4 & c_5    \\
c_6 & 0    \ebm,\]
where $c_1, c_2, \ldots, c_6$ are arbitrary nonzero real numbers.
By \eqref{eq:R}, we obtain
\[ R = CL = \bbm  c_1 c_4  & c_2 c_4+ c_3 c_5\\ c_1 c_6 & c_2 c_6\ebm.\]
It turns out that $R$ has full column rank. Indeed, the determinant of $R$ is equal to $ - c_1 c_3 c_5 c_6$ which is always nonzero.
Consequently, the FDI problem  for $(\calA,\calL,\calC)$ is solvable.
This provides a counterexample for the necessity of the condition in Theorem \ref{t:rankoutputseparable}.
\eex

 \section{A graph-theoretic condition for solvability of the FDI problem}\label{s:graph}
 
 So far, we have provided a sufficient condition for solvability of the FDI problem for $(\calA,\calL,\calC)$ in terms of the full column rank property of its associated matrix $\calR$.
 However, given such a matrix $\calR$, it is not clear how to check its full column rank property.
 Hence, in this section, we will provide a graph-theoretic condition under which a given pattern matrix $\calR$ has full column rank.
 Clearly, by Theorem \ref{t:rankoutputseparable} this will immediately lead to a graph-theoretic condition for solvability of the FDI problem for $(\calA,\calL,\calC)$.
 
We will now first review the concept of graph associated with a given pattern matrix, and the color change rule that acts on this graph.
 For more details, see \cite{JWTC2019}. 
 
 For a given pattern matrix $\mathcal{M} \in \{0,\ast,?\}^{r \times s}$ with $r \leq s$, the graph $G(\mathcal{M}) = (V,E)$ associated with $\calM$ is defined as follows.
 Take as node set $V = \{1, 2, \ldots, r\}$  and define the edge set $E \subseteq V \times V$ such that $(j,i) \in E$ if and only if $\mathcal{M}_{ij} = \ast$ or $\mathcal{M}_{ij} =?$. 
 Also, in order to distinguish between $\ast$ and $?$ entries in $\mathcal{M}$, we define two subsets $E_\ast$ and $E_?$ of the edge set $E$ as follows: $(j,i) \in E_\ast$ if and only if $\mathcal{M}_{ij} = \ast$ and $(j,i) \in E_?$ if and only if $\mathcal{M}_{ij} = ?$. 
 Then, obviously, $E = E_\ast \cup E_?$ and $E_\ast \cap E_? = \emptyset$.
 To visualize this,  solid and dashed arrows are used to represent edges in $E_\ast$ and $E_?$, respectively.  
 We say that $\calM$ has {\em full row rank} if the matrix $M$ has full row  rank for all $M  \in \calP(\calM)$.
Next, we introduce a  so-called {\em color change rule} which is defined as follows. 
 \ben
 \item Initially, color all nodes in $G(\mathcal{M})$ white.
 \item\label{i:color} If a node $i$ has exactly one white out-neighbor $j$ and $(i,j) \in E_*$,  change the color of $j$ to black.
 \item Repeat step \ref{i:color} until no more nodes can be colored black$.$
 \een
 The graph $G(\mathcal{M}) $ is called {\em colorable\/} if the nodes $1, 2, \ldots, r$ are colored black following the procedure above. Note that the remaining nodes $r+1, r+2, \ldots ,s$ can never be colored black since they have no incoming edges.
A criterion for the full row rank property of $\calM$ is then given by the following proposition.
 \bprop \cite[Theorem 11]{JWTC2019}
 	\label{t:rank}
 	Let $\mathcal{M} \in \{0,\ast,?\}^{r \times s}$ be a pattern matrix with $r \leq s$. Then $\calM$ has full row rank if and only if $G(\calM)$ is colorable.
\eprop

 Define the {\em transpose} of $\calR$ as the pattern matrix $\calR^{\top} \in \{0,\ast,?\}^{s \times r}$ with $(\calR^{\top})_{ij} = \calR_{ji}$ for $i = 1, 2, \ldots, s$ and $j = 1, 2,\ldots,r$.
 We then obtain the following obvious fact:
 \blem \label{l:columnrank}
 Consider the system $(\calA,\calL,\calC)$.  Let $\calR$ be the pattern matrix given by
 \eqref{eq:calR} and $\calR^\top$ be its transpose. Then $\calR$ has full column rank if and only if $G(\calR^\top)$ is colorable.
 \elem

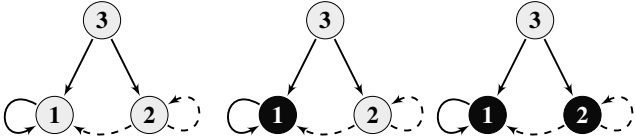
\begin{figure}[h!]
 	\centering
 	\begin{subfigure}{0.16\textwidth}
 		\centering
 		\scalebox{0.9}{\begin{tikzpicture}[scale=0.35]
 		\tikzset{VertexStyle1/.style = {shape = circle,
 				color=black,
 				fill=white!93!black,
 				minimum size=0.5cm,
 				text = black,
 				inner sep = 2pt,
 				outer sep = 1pt,
 				minimum size = 0.55cm}
 		}	
 		\node[VertexStyle1,draw](1) at (-2,0) {$\bf 1$};
 		\node[VertexStyle1,draw](2) at (2,0) {$\bf 2$};
 		\node[VertexStyle1,draw](3) at (0,4) {$\bf 3$};
 		\Loop[style={> = latex',->, out=150, in=-150,line width=0.8pt}, dist=2cm](1)
 		\Loop[style={> = latex',->, out=-30, in=30,line width=0.8pt,dashed}, dist=2cm](2)
 		\Edge[style = {->,> = latex',pos = 0.7},color=black
 		, labelstyle={inner sep=0pt}](3)(2);
 		\Edge[style = {->,> = latex',pos = 0.7},color=black
 		, labelstyle={inner sep=0pt}](3)(1);
 		\Edge[style = {->,> = latex',pos = 0.7, out = -150, in = -30, dashed},color=black
 		, labelstyle={inner sep=0pt}](2)(1);
 		\end{tikzpicture}}
 		\caption{The initial graph.}
 		\label{g:R1}
 	\end{subfigure}
 	\begin{subfigure}{0.15\textwidth}
 		\centering
 		\scalebox{0.9}{\begin{tikzpicture}[scale=0.35]
 		\tikzset{VertexStyle1/.style = {shape = circle,
 				color=black,
 				fill=white!93!black,
 				minimum size=0.5cm,
 				text = black,
 				inner sep = 2pt,
 				outer sep = 1pt,
 				minimum size = 0.55cm}
 		}	
 		\tikzset{VertexStyle2/.style = {shape = circle,
 				color=black,
 				fill=black!96!white,
 				minimum size=0.5cm,
 				text = white,
 				inner sep = 2pt,
 				outer sep = 1pt,
 				minimum size = 0.55cm}
 		}
 		\node[VertexStyle2,draw](1) at (-2,0) {$\bf 1$};
 	\node[VertexStyle1,draw](2) at (2,0) {$\bf 2$};
 	\node[VertexStyle1,draw](3) at (0,4) {$\bf 3$};
 		\Loop[style={> = latex',->, out=150, in=-150,line width=0.8pt}, dist=2cm](1)
 		\Loop[style={> = latex',->, out=-30, in=30,line width=0.8pt,dashed}, dist=2cm](2)
 		\Edge[style = {->,> = latex',pos = 0.7},color=black
 		, labelstyle={inner sep=0pt}](3)(2);
 		\Edge[style = {->,> = latex',pos = 0.7},color=black
 		, labelstyle={inner sep=0pt}](3)(1);
 		\Edge[style = {->,> = latex',pos = 0.7, out = -150, in = -30, dashed},color=black
 		, labelstyle={inner sep=0pt}](2)(1);
 		\end{tikzpicture}}
 		\caption{ Node $1$ colors $1$.}
 		\label{g:R2}
 	\end{subfigure}
 \begin{subfigure}{0.15\textwidth}
 		\centering
 		\scalebox{0.9}{\begin{tikzpicture}[scale=0.35]
 		\tikzset{VertexStyle1/.style = {shape = circle,
 				color=black,
 				fill=white!93!black,
 				minimum size=0.5cm,
 				text = black,
 				inner sep = 2pt,
 				outer sep = 1pt,
 				minimum size = 0.55cm}
 		}	
 		\tikzset{VertexStyle2/.style = {shape = circle,
 				color=black,
 				fill=black!96!white,
 				minimum size=0.5cm,
 				text = white,
 				inner sep = 2pt,
 				outer sep = 1pt,
 				minimum size = 0.55cm}
 		}
 		\node[VertexStyle2,draw](1) at (-2,0) {$\bf 1$};
 	\node[VertexStyle2,draw](2) at (2,0) {$\bf 2$};
 	\node[VertexStyle1,draw](3) at (0,4) {$\bf 3$};
 		\Loop[style={> = latex',->, out=150, in=-150,line width=0.8pt}, dist=2cm](1)
 		\Loop[style={> = latex',->, out=-30, in=30,line width=0.8pt,dashed}, dist=2cm](2)
 		\Edge[style = {->,> = latex',pos = 0.7},color=black
 		, labelstyle={inner sep=0pt}](3)(2);
 		\Edge[style = {->,> = latex',pos = 0.7},color=black
 		, labelstyle={inner sep=0pt}](3)(1);
 		\Edge[style = {->,> = latex',pos = 0.7, out = -150, in = -30, dashed},color=black
 		, labelstyle={inner sep=0pt}](2)(1);
 		\end{tikzpicture}}
 		\caption{ Node $3$ colors $2$.}
 		\label{g:R3}
 	\end{subfigure}
 	\caption{The graph $G(\calR^\top)$ is colorable.}
 	\label{g:R}
 \end{figure}

This then immediately yields the main result of  this section which provides a graph-theoretic condition under which the FDI problem for  $(\calA,\calL,\calC)$ is solvable.
 \bthm \label{t:graph}
 Consider the system $(\calA,\calL,\calC)$.  Suppose that the indices $\eta_i$ exists for $i = 1,2, \ldots, q$. Let $\calR$ be the pattern matrix given by
 \eqref{eq:calR}.
 Then, the FDI problem for  $(\calA, \calL, \calC)$ is solvable if $G(\calR^\top)$ is colorable.
 \ethm
 \BP
 The proof follows immediately from Theorem \ref{t:rankoutputseparable} and Lemma \ref{l:columnrank}.
 \EP
 To conclude this section, we will  provide an example.
 \addtolength{\textheight}{-.67cm}
 \bex \label{ex:4} 
 Consider the system $(\calA, \calL, \calC)$ with
 \[ \calA = \resizebox{0.27\hsize}{!}{ $\bbm \ast & 0   & 0    & 0   & 0\\
 \ast & ?   & 0    & ?   & 0\\ 
 0    &\ast & \ast & ?   & 0\\ 
 \ast &0    &0     & ?   &\ast\\
 0    &0    &\ast  & 0   &\ast\ebm$},~~ 
 \calL = \resizebox{0.11\hsize}{!}{ $\bbm \ast & 0  \\ 
 ? &\ast\\
 0 & 0  \\
 0 & 0 \\
 0 & 0 \ebm$},~~ \calC = \resizebox{0.27\hsize}{!}{ $\bbm 0    & 0   & 0    & \ast & 0\\
 0    & 0   & 0    & ?    & ?\\ 
 0    & 0   & 0    & \ast & \ast \ebm.$} \]
 By multiplying the pattern matrices, we obtain that
 \[\bbm \calC \calL_1 & \calC \calA \calL_1\ebm = \resizebox{0.12\hsize}{!}{ $\bbm 0 & \ast   \\ 
 0 & ? \\
 0 & \ast\ebm$}\]
 and                             
 \[\bbm \calC \calL_2 & \calC \calA \calL_2 & \calC \calA^2 \calL_2\ebm 
 =  \resizebox{0.18\hsize}{!}{ $\bbm 0 & 0 & 0\\ 
 0 & 0 & ?\\
 0 & 0 &\ast\ebm$},\] where $\calL_i$ is the $i$th column of $\calL$.
 By \eqref{eq:calR}, it follows that the associated matrix $\calR$ and its transpose $\calR^\top$ are given by 
 $$\calR = \bbm \calC \calA \calL_1&\calC \calA^2 \calL_2 \ebm 
 = \resizebox{0.12\hsize}{!}{ $\bbm  \ast & 0 \\
 ? & ?\\
 \ast & \ast 
 \ebm$}$$
 and 
 $$\calR^\top= \resizebox{0.18\hsize}{!}{ $\bbm  \ast & ? & \ast\\
 0 & ? & \ast
 \ebm$}.$$
 As depicted  in Fig. \ref{g:R}  $G(\calR^\top)$ is colorable.
 Indeed, initially let all nodes in $G(\calR^\top)$ be colored white as shown in Fig. \ref{g:R}(a). Node $1$ then colors itself black as depicted in  Fig. \ref{g:R}(b), and finally node $3$ colors $2$ to black as in  Fig. \ref{g:R}(c).
 Therefore, by Theorem \ref{t:graph}, the FDI problem for  $(\calA,\calL,\calC)$ is solvable.
\eex

\section{Conclusion} \label{s:conclusions}
In this paper, we have studied the FDI problem for linear structured systems. 
We have established a necessary and sufficient condition for solvability of the FDI problem for a given particular LTI system.
Based on this, we have established a necessary condition under which the FDI problem for structured systems is solvable.
Moreover, we have developed a sufficient condition for solvability of the FDI problem in terms of a rank test on a pattern matrix associated with the structured system.
Next, we have provided a counterexample to show that this condition is not necessary.
Finally, we have developed a graph-theoretic condition for solvability of the FDI problem using the concept of colorability of a graph.

This paper has only established  a necessary condition and sufficient conditions for solvability of the FDI problem for structured systems.
Finding necessary {\em and} sufficient conditions for solvability of the FDI problem is still an open problem.
In addition, as we have mentioned in section \ref{s:GFDI}, this paper does not consider the stability of the unknown input observers.
Therefore, another possible future research direction is to establish conditions under which stable unknown input observers exist for linear structured systems.
Furthermore, investigating solvability of FDI for structured systems with constraints, such as allowing dependencies on some nonzero and arbitrary entries \cite{LM2019, JTWC2018, JTBC20181}, is also a possibility for future research.

\bibliographystyle{IEEEtran}
\bibliography{CDC2020}

\end{document}